\documentclass[a4paper]{amsart}

\usepackage{amsfonts}
\usepackage{amsmath}
\usepackage{amssymb}
\usepackage{mathrsfs}
\newtheorem {Def}{Definition}[section]
\newtheorem {Lem}[Def]{Lemma}
\newtheorem {Thm}[Def]{Theorem}
\newtheorem {Prop}[Def]{Proposition}
\newtheorem {Cor}[Def]{Corollary}
\newcommand{\CC}{{\mathbb C}}
\newcommand{\id}{\mathrm{id}}
\newcommand{\comp}{\!\circ\!}
\newcommand{\bez}{\setminus}
\newcommand{\eps}{\varepsilon}
\newcommand{\ph}{\varphi}
\newcommand{\tens}{\otimes}
\newcommand{\Sum}{\sum\limits}
\newcommand{\M}[1]{M\!\left(#1\right)}
\newcommand{\Span}[1]{\mathrm{span}\left\{#1\right\}}
\newcommand{\Ci}[1]{C_{\infty}\!\left(#1\right)}
\newcommand{\refeq}[1]{{\rm (\ref{#1})}}
\newcommand{\Aa}{{\mathscr A}}
\newcommand{\Bb}{{\mathscr B}}
\newcommand{\Aah}{{\Aa\sp{\:\,\!\!\widehat{}}}}
\newcommand{\Aahz}{{\Aah}^{\:0}}
\newcommand{\Ap}{\mathscr{A\!P}}
\newcommand{\sh}{{\,\sharp}}
\newcommand{\st}{{\,*}}
\newcommand{\II}{{\mathbb I}}
\newcommand{\Ff}{{\mathscr F}}
\newcommand{\del}{\delta}
\newcommand{\kap}{\kappa}
\newcommand{\delhat}{{\del\sp{\!\widehat{}}}\,}

\begin{document}

\title{Compactifications of Discrete Quantum Groups}
\author{Piotr Miko\l{}aj So\l{}tan}
\email{Piotr.Soltan@fuw.edu.pl}\thanks{
Partially supported by Komitet Bada\'n Naukowych grant no.~2P03A04022, the
Foundation for Polish Science and Deutsche Forschungsgemeinschaft.}
\address{Department of Mathematical Methods in Physics,\\
Faculty of Physics, University of Warsaw}
\begin{abstract}
Given a discrete quantum group $(\Aa,\del)$ we construct Hopf
$*$-algebra $\Ap$ which is a unital $*$-subalgebra of the multiplier algebra of
$\Aa$. The structure maps for $\Ap$ are inherited from $\M{\Aa}$ and thus the
construction yields a compactification of $(\Aa,\del)$ which is analogous to
the Bohr compactification of a locally compact group. This algebra has the
expected universal property with respect to homomorphisms from multiplier Hopf
algebras of compact type (and is therefore unique). This provides an easy proof
of the fact that for a discrete quantum group with an infinite dimensional
algebra the multiplier algebra is never a Hopf algebra.
\end{abstract}
\maketitle

\section{Introduction}

The research presented in this paper was motivated by investigations in the
theory of ${\mathrm C}^*$-algebraic quantum groups (\cite{kv,mnw}). This theory
is
based on the theory of locally compact groups. Let $G$ be a locally compact
group and let $A=\Ci{G}$ the algebra of all continuous functions vanishing at
infinity on $G$. Then the algebra of bounded continuous functions on $G$ is
naturally isomorphic to the multiplier algebra $\M{A}$ of $A$. A function
$f\in\M{A}$ is almost periodic on $G$ if and only if
$\del(f)\in\M{A}\tens\M{A}$, where $\del$ is a morphism from $A$ to
$A\tens A=\Ci{G\times G}$ given by $\bigl(\del{f}\bigr)(g_1,g_2)=f(g_1g_2)$
and ``$\tens$'' stands for the completed tensor product of ${\mathrm
C}^*$-algebras. The
set of almost periodic functions is a commutative unital
${\mathrm C}^*$-subalgebra of $\M{A}$. This algebra is the algebra of continuous
functions on the Bohr compactification of $G$ (cf.~\cite[\S 41]{loom}). Our aim
is to generalize this construction to discrete quantum groups in the framework
of multiplier Hopf algebras introduced by Van Daele in \cite{mha}.

We shall now describe the contents of the paper. In section \ref{prel} we
gather necessary information about multiplier Hopf algebras and discrete
quantum groups. Section \ref{sli} deals with the concept of slices with
reduced functionals which is used in the next section. In Section \ref{cdqg}
we define the algebra of almost periodic elements for a discrete quantum group
and show that it is a Hopf $*$-algebra. We also prove that this Hopf algebra
has a universal property for morphisms of multiplier Hopf algebras of compact
type into the original discrete quantum group. At the end of Section
\ref{cdqg} we include some corollaries of our construction.

\section{Preliminaries}\label{prel}

All algebras we shall consider in this paper will be over the field of complex
numbers. Let $\Aa$ be an algebra with non degenerate product. By $\M{\Aa}$ we
shall denote the multiplier algebra of $\Aa$ (see \cite{mha}). We shall use the
fact that a multiplier $m$ is determined by the linear map corresponding to
multiplication by $m$ from the left. Indeed: if $m,n\in\M{\Aa}$ and $ma=na$ for
all $a\in\Aa$ then $(bm)a=b(ma)=b(na)=(bn)a$ for any $b$. By non degeneracy of
the product in $\Aa$ we have that right multiplication by $n$ and $m$ also
agree as linear maps on $\Aa$. In other words, $m=n$ as multipliers.

Let $\Aa$ and $\Bb$ be algebras with non degenerate products. A homomorphism
$\Phi\colon\Aa\to\M{\Bb}$ is non degenerate if $\Phi(\Aa)\Bb=\Bb\Phi(\Aa)=\Bb$.
Such a non degenerate homomorphism has a unique extension to a homomorphism of
unital algebras $\M{\Aa}\to\M{\Bb}$. In particular non degenerate
homomorphisms can be composed just as usual homomorphisms between algebras. A
composition of non degenerate homomorphisms is non degenerate.

We shall be concerned with the theory of multiplier Hopf algebras developed in
\cite{mha,afgd}. Recall that a multiplier Hopf algebra is a pair $(\Aa,\del)$
consisting of an algebra $\Aa$ with non degenerate product and a homomorphism
$\del\colon\Aa\to\M{\Aa\tens\Aa}$ such that the maps
\[
\begin{array}{r@{\;\colon\Aa\tens\Aa\ni(a\tens b)\longmapsto\;}l@{\smallskip}}
T_1&\del(a)(I\tens b),\\
T_2&(a\tens I)\del(b)
\end{array}
\]
have range equal to $\Aa\tens\Aa$ and are bijective and such that the linear
maps on $\Aa\tens\Aa\tens\Aa$ given by $T_2\tens\id$ and $\id\tens T_1$
commute (in other words, $\del$ is {\em coassociative}\/). It follows from this
definition that $\del$ is a non degenerate homomorphism.

In the fundamental reference \cite{mha} it is shown that for a multiplier Hopf
algebra $(\Aa,\del)$ there exists a multiplicative functional $\eps$ on $\Aa$
called {\em counit}\/ such that for all $a,b\in\Aa$
\[
(\id\tens\eps)\bigl((a\tens I)\del(b)\bigr)=ab
=(\eps\tens\id)\bigl((\del(a)(I\tens b)\bigr).
\]
This combined with bijectivity of either $T_1$ or $T_2$ implies that
$\Aa^2=\Aa$ which proves to be a very useful fact. Moreover
$(\id\tens\del)$ and $(\del\tens\id)$ are non degenerate homomorphisms
$\Aa\tens\Aa\to\M{\Aa\tens\Aa\tens\Aa}$. In particular we can consider the
compositions $(\del\tens\id)\comp\del$ and $(\id\tens\del)\comp\del$ form
$\Aa$ to $\M{\Aa\tens\Aa\tens\Aa}$ and their extensions
$\M{\Aa}\to\M{\Aa\tens\Aa\tens\Aa}$. In this context the coassociativity means
that the two latter maps are equal:
\[
(\id\tens\del)\del=(\del\tens\id)\del.
\]

There also exists an antihomomorphism $\kap\colon\Aa\to\M{\Aa}$ called {\em
coinverse}\/ (or {\em antipode}\/) such that for all $a,b,c\in\Aa$
\[
\begin{array}{r@{\;=\;}l@{\smallskip}}
m\left[(\id\tens\kap)\bigl((a\tens I)\del(b)\bigr)(I\tens c)\right]&a\eps(b)c,\\
m\left[(a\tens I)(\kap\tens\id)\bigl(\del(b)(I\tens c)\bigr)\right]&a\eps(b)c.
\end{array}
\]

If $(\Aa,\delta)$ is a multiplier Hopf algebra and in addition $\Aa$ is a
$*$-algebra and $\del$ is a $*$-homomorphism, the pair $(\Aa,\del)$ is then
called a multiplier Hopf $*$-algebra. The counit is then a $*$-character of
$\Aa$ and the coinverse satisfies $S\bigl(S(a)^*\bigr)^*=a$ for all $a\in\Aa$.

For a multiplier Hopf algebra $(\Aa,\del)$ we can consider the pair
$(\Aa,\del')$ where $\del'$ is a composition of $\del$ with an extension to
$\M{\Aa\tens\Aa}$ of the flip map $a\tens b\mapsto b\tens a$. If
$(\Aa,\del')$ is also a multiplier Hopf algebra then $(\Aa,\del)$ is called
{\em regular}. The coinverse of a regular multiplier Hopf algebra is
necessarily an antiisomorphism $A\to A$. In particular it is non degenerate
(as a homomorphism from $\Aa$ to the opposite algebra of $\M{\Aa}$) and
extends to an antiisomorphism $\M{\Aa}\to\M{\Aa}$.

A regular multiplier Hopf algebra $(\Aa,\del)$ is said to be of {\em discrete
type}\/ if there exists a non zero element $h\in\Aa$ such that $ah=\eps(a)h$
for all $a\in\Aa$ (cf.~\cite{dt}). A  typical example of a regular multiplier
Hopf algebra of discrete type is a {\em discrete quantum group}, i.e.~a
multiplier Hopf $*$-algebra $(\Aa,\del)$ such that $\Aa$ is a direct sum of full
matrix algebras (\cite[Def.~2.3]{dqg}). It is well known (\cite{mha}) that
multiplier Hopf $*$-algebras are regular. Involutive structure will not play
an important role in our considerations. Regularity will be much more
important, but we shall keep the $*$-structure in order to be able to use
results from the theory of Hopf $*$-algebras (\cite{dp}).

Discrete quantum groups appeared first in \cite[Sect.~3]{pw}. They were
defined and studied in the framework of multiplier Hopf algebras in
\cite{dqg}.

A functional $\ph$ on a multiplier Hopf algebra $(\Aa,\del)$ is called {\em left
invariant}\/ if $(\id\tens\ph)\del(a)=\ph(a)I$ for all $a\in\Aa$
($(\id\tens\ph)\del$ is a map from $\Aa$ to $\M{\Aa}$). Similarly a functional
$\psi$ on $\Aa$ is {\em right invariant}\/ if $(\psi\tens\id)\del(a)=\ph(a)I$
for
all $a\in\Aa$. If $\ph$ is a left invariant functional for a regular multiplier
Hopf algebra $(\Aa,\del)$ then $\psi=\ph\comp S$ is right invariant. The theory
of regular multiplier Hopf algebras with invariant functionals is very rich.
In particular the dual regular multiplier Hopf algebra can be defined as
\[
\Aah=\bigl\{a\ph:a\in\Aa\bigr\},
\]
where $\ph$ is any non trivial left invariant functional. This definition does
not depend on the choice of $\ph$ (which is in fact unique up to rescaling)
and $\Aah$ can be endowed with a comultiplication $\delhat$ dual to
multiplication in $\Aa$ and it becomes a regular multiplier Hopf
algebra with invariant functionals. The biduality theorem
\cite[Thm.~4.12]{afgd} says that the dual of $(\Aah,\delhat)$ is naturally
isomorphic to $(\Aa,\del)$. The same assertions are true if we consider
multiplier Hopf $*$-algebras.

Let us conclude with a statement that all regular multiplier Hopf algebras of
discrete type have non trivial invariant functionals (\cite{dt}, see also
\cite{dqg}).

\section{Slices}\label{sli}

Let $\Aa$ be an algebra with non degenerate product. The space of all linear
functionals on $\Aa$ will be denoted by $\Aa^\sh$. This vector space carries a
natural $\Aa$-bimodule structure: for $f\in\Aa^\sh$ and $a\in\Aa$
\[
\begin{array}{r@{\;=\;}l@{\smallskip}}
\bigl(af\bigr)(b)&f(ba),\\
\bigl(fa\bigr)(b)&f(ab)
\end{array}
\]
for all $b\in\Aa$. The space of all {\em reduced linear functionals}\/ on $\Aa$
is by definition
\[
\Aa^\st=\Span{\bigl.afb:\:f\in\Aa^\sh,\:a,b\in\Aa\bigr.}.
\]

Now assume that $\Aa$ is a direct sum of matrix algebras. Any reduced functional
on $\Aa$ has a natural extension from $\Aa$ to $\M{\Aa}$:
\[
\bigl(afb\bigr)(m)=f(bma)
\]
for any multiplier $m$. To see that this is well defined we must show that
if $a_1,\ldots,a_n,b_1,\ldots,b_n\in\Aa$ and $f_1,\ldots,f_n\in\Aa^\sh$
then
\[
\sum f_i(a_i c b_i)=0
\]
for all $c\in\Aa$ implies that
\[
\sum f_i(a_i m b_i)=0
\]
for all $m\in\M{\Aa}$. Let $e$ be the unit of the direct sum of the
matrix algebras containing the elements $a_1,\ldots,a_n$. Then $e\in\Aa$ and
$a_ie=a_i$ for $1\leq i\leq n$. Now we see that for any $m\in\M{\Aa}$
\[
\sum f_i(a_i m b_i)=\sum f_i\bigl(a_i(em) b_i\bigr)=0,
\]
as $em\in\Aa$.

If $\Bb$ is an algebra with non degenerate product then $\Bb\tens\Aa$ is
also an algebra with non degenerate product (\cite[Lemma A.2]{mha}). For any
$\xi\in\Aa^\st$ and $\zeta\in\Bb^\st$ the tensor product $\zeta\tens\xi$ is a
reduced functional on $\Bb\tens\Aa$ and, as such, extends to $\M{\Bb\tens\Aa}$.

\begin{Prop}
Let $\Aa$ and $\Bb$ be algebras with non degenerate products and let $Y$ be a
multiplier of $\Bb\tens\Aa$. Then for any $\xi\in\Aa^\st$ there exists a unique
multiplier $m\in\M{\Bb}$ such that
\begin{equation}\label{slice}
(\zeta\tens\xi)(Y)=\zeta(m)
\end{equation}
for all $\zeta\in\Bb^\st$. The multiplier $m$ is called a\/ {\em right slice
of $Y$ with $\xi$} and will be denoted by $(\id\tens\xi)(Y)$.
\end{Prop}

\noindent\begin{proof}
Let $\xi=afc$ with $a,c\in\Aa$ and $f\in\Aa^\sh$. Define left and
right multiplication by $m$ as
\[
\begin{array}{r@{\;=\;}l@{\smallskip}}
m\,b_2&(\id\tens f)\bigl((I\tens c)Y(b_2\tens a)\bigr),\\
b_1m&(\id\tens f)\bigl((b_1\tens c)Y(I\tens a)\bigr).
\end{array}
\]
for $b_1,b_2\in\Bb$. It remains to prove that
\[
(b_1m)b_2=b_1(m\,b_2)
\]
and this follows from associativity of multiplication in $\Bb\tens\Aa$:
\[
\begin{array}{r@{\;=\;}l@{\smallskip}}
b_1\left[(\id\tens f)\bigl((I\tens c)Y(b_2\tens a)\bigr)\right]
&(\id\tens f)\bigl((b_1\tens I)(I\tens c)Y(b_2\tens a)\bigr)\\
&(\id\tens f)\bigl((b_1\tens c)Y(b_2\tens a)\bigr)\\
&(\id\tens f)\bigl((b_1\tens c)Y(I\tens a)(b_2\tens I)\bigr)\\
&\left[(\id\tens f)\bigl((b_1\tens c)Y(I\tens a)\bigr)\right]b_2.
\end{array}
\]
Formula \refeq{slice} holds, for if $\zeta=b_2f'b_1$ with $b_1,b_2\in\Bb$ and
$f'\in\Bb^\sh$ then
\[
\begin{array}{r@{\;=\;}l@{\smallskip}}
\zeta\bigl((\id\tens\xi)(Y)\bigr)&f'\bigl(b_1(\id\tens\xi)(Y)b_2\bigr)\\
&f'\left[(\id\tens f)\bigl((b_1\tens c)Y(b_2\tens a)\bigr)\right]\\
&(f'\tens f)\bigl((b_1\tens c)Y(b_2\tens a)\bigr)\\
&(b_2f'b_1\tens afc)(Y).
\end{array}
\]

Let $n$ be another multiplier of $\Bb$ such that for any $\zeta\in\Bb^\st$
\[
(\zeta)(n)=(\zeta\tens\xi)(Y).
\]
Then for any $g\in\Bb^\sh$ and all $b,b'\in\Bb$ we have $g(bnb')=g(bmb')$ and it
follows that $bnb'=bmb'$ for all $b,b'\in\Bb$. Therefore $b(nb'-mb')=0$ for all
$b\in\Bb$. Since $nb'-mb'\in\Bb$, this implies that $nb'=mb'$ for all $b'$ and
this means that $m=n$ as linear maps on $\Bb$.
\end{proof}

\section{Compactification of a discrete quantum group}\label{cdqg}

The algebra of $n\times n$ matrices with complex entries will be denoted by
$M_n$. Let $\Aa$ be a direct sum of a family of such full matrix algebras. Then
$\M{\Aa}$ is a product if the same family of matrix algebras. To see the
isomorphism let $\Aa=\bigoplus\limits_{\iota\in\II}\Aa_\iota$ with
$\Aa_\iota=M_{n_\iota}$ for each $\iota\in\II$. The action of an infinite family
$(m_\iota)_{\iota\in\II}$ with $m_\iota\in M_{n_\iota}$ on elements of $\Aa$ is,
of course, given by matrix multiplication in each summand. Conversely for a
multiplier $m$ of $\Aa$ the corresponding family $(m_\iota)_{\iota\in\II}$ is
obtained by setting $m_\iota=me_{\iota}$ where $e_\iota$ is the unit of the
matrix algebra $\Aa_\iota$.

\begin{Lem}\label{L}\sloppy
Let $\Aa$ be a direct sum of matrix algebras and let $x_1,\ldots,x_N$ be
linearly independent multipliers of $\Aa$. Then there exists an element
$e\in\Aa$ such that $x_1e,\ldots,x_Ne$ are linearly independent elements of
$\Aa$. The element $e$ may be chosen to be a central idempotent.
\end{Lem}

\noindent\begin{proof}
To proceed we must introduce notation
\[
\Aa=\bigoplus_{\iota\in\II}\Aa_\iota
\]
with $\Aa_\iota=M_{n_\iota}$.

\sloppy
Suppose, contrary to the statement of the lemma, that for any $a\in\Aa$ the set
$\{x_1c,\ldots,x_Nc\}$ is not linearly independent. Let $\Ff$ be the family of
finite subsets the index set $\II$. Notice that the family $\Ff$ is directed by
inclusion. For any $F\in\Ff$ let $e_F$ be the unit of the finite dimensional
algebra $\bigoplus\limits_{\iota\in F}\Aa_\iota$. We are under assumption that
\begin{equation}\label{assumpt}
\left(\begin{array}{c@{\smallskip}}
\text{For any $F\in\Ff$ there exists a vector}\\
(\lambda_1^F,\ldots,\lambda_N^F)\in\CC^N\bez\{0\}\\
\text{such that }\Sum_{k=1}^{N}\lambda_k^Fx_ke_F=0.
\end{array}\right)
\end{equation}
For any $F\in\Ff$ we shall denote by $V_F$ the (non zero) subspace of $\CC^N$
consisting of all vectors $(\lambda_1^F,\ldots,\lambda_N^F)$ fulfilling the
formula in \refeq{assumpt}.

Let $F_0,F\in\Ff$ with $F_0\subset F$. Then for any
$(\mu_1,\ldots,\mu_N)\in V_F$ we have
\[
\sum_{k=1}^{N}\mu_kx_ke_F=0
\]
and multiplying this relation from the right by $e_{F_0}$ we obtain
\[
\Sum_{k=1}^{N}\mu_kx_ke_{F_0}=0.
\]
This means that $(\mu_1,\ldots,\mu_N)\in V_{F_0}$ and consequently
$V_F\subset V_{F_0}$. Let
$V_{\infty}=\bigcap\limits_{F\in\Ff}V_F$. A moment of reflection shows that
this subspace is non zero.

Let $(\alpha_1,\ldots,\alpha_N)$ be a non zero vector in $V_\infty$. then for
any $F\in\Ff$ we have
\[
\sum_{k=1}^{N}\alpha_kx_ke_F=0.
\]
Now for any $a\in\Aa$ there is an $F\in\Ff$ such that $ae_F=e_Fa=a$. It
follows that for any $a\in\Aa$
\[
\sum_{k=1}^{N}\alpha_kx_ka=\sum_{k=1}^{N}\alpha_kx_kae_F
=\sum_{k=1}^{N}\alpha_kx_ke_Fa=\left(\sum_{k=1}^{N}\alpha_kx_ke_F\right)a=0.
\]
In other words $\Sum_{k=1}^{N}\alpha_kx_k=0$, i.e.~the elements
$\{x_1,\ldots,x_N\}$ are not linearly independent in $\M{\Aa}$.

This contradiction means that \refeq{assumpt} is not true. Consequently there
exists a central projection $e\in\Aa$ such that $x_1e,\ldots,x_Ne$ are linearly
independent.
\end{proof}

\begin{Def}\label{DefAP}
Let $(\Aa,\del)$ be a discrete quantum group. The set $\Ap$ of almost
periodic elements for $(\Aa,\del)$ is defined as
\[
\Ap=\bigl\{x\in\M{\Aa}:\:\del(x)\in\M{\Aa}\tens\M{\Aa}\bigr\}.
\]
\end{Def}

\begin{Thm}
Let $(\Aa,\del)$ be a discrete quantum group and denote its coinverse by $\kap$.
Then $\Ap$ defined in Definition \ref{DefAP} is a unital $*$-subalgebra of
$\M{\Aa}$. Moreover we have
\begin{description}
\item[1.] $\del(\Ap)\subset\Ap\tens\Ap$,
\item[2.] $\kap(\Ap)=\Ap$.
\end{description}
\end{Thm}

\noindent\begin{proof}
The set $\Ap$ is a unital $*$-subalgebra of $\M{\Aa}$ as a pre image of a
unital $*$-subalgebra in a $*$-homomorphism of unital algebras. Let us
concentrate on the comultiplication. Fix $x\in\Ap$. Then we can write
\begin{equation}\label{suma1}
\del(x)=\sum_{k=1}^{N}x_{k}\tens y_k
\end{equation}
with $(x_k)_{k=1,\ldots,N}$ and $(y_k)_{k=1,\ldots,N}$ in $\M{\Aa}$ and
$(y_k)_{k=1,\ldots,N}$ linearly independent. By coassociativity of $\del$ we
have
\begin{equation}\label{coas}
\sum_{k=1}^{N}\del(x_k)\tens y_k=\sum_{k=1}^{N}x_{k}\tens\del(y_k)
\end{equation}
It follows that
$\del(x)\in\M{\Aa\tens\Aa}\tens\M{\Aa}\cap\M{\Aa}\tens\M{\Aa\tens\Aa}$. For
any reduced functional $\xi$ on $\Aa$ we can apply the map
\begin{equation}\label{form}
\id_{\Aa\tens\Aa}\tens\xi=\id_{\Aa}\tens(\id\tens\xi)
\end{equation}
to $\del(x)$. By \refeq{coas} and \refeq{form} the value will lie in
$\M{\Aa}\tens\M{\Aa}$.

Let $e$ be a central idempotent in $\Aa$ such that $y_1e,\ldots,y_Ne$ are
linearly independent (as in Lemma \ref{L}). Since the ideal generated by $e$ is
a finite dimensional direct summand of $\Aa$, it is easy to define functionals
$\xi_1,\ldots,\xi_N$ on
$\Aa$ such that
\[
\xi_l(y_ke)=\delta_{k,l}.
\]
for $k,l=1,\ldots,N$ and
\[
\xi_k=e\xi_k
\]
for $k=1,\ldots,N$. Thus the functionals $(\xi_k)_{k=1,\ldots,N}$ are reduced
and taking slices of \refeq{coas} with $\xi_l$ we obtain
\[
\begin{array}{l@{\smallskip}}
\M{\Aa}\tens\M{\Aa}\ni(\id_{\Aa\tens\Aa}\tens\xi_l)
\left(\Sum_{k=1}^{N}\del(x_k)\tens y_k\right)\\
\qquad\qquad\qquad\qquad=\Sum_{k=1}^{N}\del(x_k)\xi_l(y_k)=\del(x_l).
\end{array}
\]
This means that for each $l\in\{1,\ldots,N\}$ we have $x_l\in\Ap$. In
particular $\del(x)\in\Ap\tens\M{\Aa}$. Now choosing maximal linearly
independent subset out of the elements $\{x_1,\ldots,x_N\}$ we can rearrange
the sum \refeq{suma1} to have linearly independent elements of $\Ap$ making
up the left leg of $\del(x)$. Then we can show that all elements making up
the right leg are contained in $\Ap$ using the same technique as we have done
for the $\{x_1,\ldots,x_N\}$. The final statement being
\[
\del(x)\in\Ap\tens\Ap,
\]
which proves \textbf{1.}

In what follows we shall use techniques which have become standard in the
theory of multiplier Hopf algebras. The relation between coinverse and
comultiplication in a multiplier Hopf algebra is given in
\cite[Prop.~5.6]{mha}: for any $a,b\in\Aa$ we have
\begin{equation}\label{koinv}
\bigl(I\tens\kap(b)\bigr)\del\bigl(\kap(a)\bigr)=(\kap\tens\kap)
\bigl(\del'(a)(I\tens b)\bigr)
\end{equation}
and $\del'$ is the composition of $\del$ with the flip automorphism
(cf.~Section \ref{prel}).

Now our $(\Aa,\del)$ is a regular multiplier Hopf algebra. In particular
$\kap\tens\kap$ is an antiisomorphism of $\Aa\tens\Aa$ onto itself and it
extends to an antiisomorphism of $\M{\Aa\tens\Aa}$ onto itself. Thus the right
hand side of \refeq{koinv} reads
\[
\bigl(I\tens\kap(b)\bigr)(\kap\tens\kap)\del'(a).
\]
In view of non degeneracy of the product we conclude that
\[
\del\bigl(\kap(a)\bigr)=(\kap\tens\kap)\del'(a)
\]
for all $a\in\Aa$. Thus the two non degenerate antihomomorphisms
$\del\comp\kap$ and $(\kap\tens\kap)\comp\del'$ agree on $\Aa$ and
their (unique) extensions agree on $\M{\Aa}$. With this identity point
\textbf{2.} is trivial.
\end{proof}

In order to fully describe the rest of the structure of $\Ap$ we shall use
the fact that the dual multiplier Hopf algebra of $(\Aa,\del)$ is a Hopf
$*$-algebra. First let us observe that $\M{\Aa}$ can be identified with a
subspace of $\Aah^\sh$. Indeed: any element of $\Aah$ is of the form $a\ph$
for some $a\in\Aa$ (with $\ph$ a left invariant functional on $\Aa$).
Therefore any $m\in\M{\Aa}$ determines a linear map
\[
\Aah\ni a\ph\longmapsto\ph(ma)\in\CC.
\]
Using the biduality theorem (\cite[Thm.~4.12]{afgd}) one can see that
the structure of a $*$-algebra on $\M{\Aa}$ coincides with the one obtained
from $\Aah$ on $\Aah^\sh$ as described in the beginning of \cite[Sect.~3]{dp}.
Now it is easy to see that $\Ap$ is a subspace of the space $\Aahz$ defined
in \cite[Def.~3.2]{dp}. Moreover using the biduality theorem one can see that
the maps $\del$, $\eps$ and $\kap$ on $\Ap$ are exactly the restrictions of
the corresponding maps making $\Aahz$ a Hopf $*$-algebra (as in
\cite[Thm.~3.3]{dp}). In particular we have
\[
\begin{array}{r@{\;=\;}l@{\smallskip}}
m\bigl((\kap\tens\id)\del(x)\bigr)&\eps(x)I,\\
m\bigl((\id\tens\kap)\del(x)\bigr)&\eps(x)I
\end{array}
\]
for all $x\in\Ap$.

We can summarize the above considerations in the following theorem:

\begin{Thm}\label{main}
Let $(\Aa,\del)$ be a discrete quantum group. Then the algebra $\Ap$ defined
in Definition \ref{DefAP} with structure maps inherited from $\M{\Aa}$ is a
Hopf $*$-algebra.
\end{Thm}

$\Ap$ is a unital $*$-subalgebra of the multiplier algebra of $\Aa$. This means,
in particular, that the inclusion
$\Ap\hookrightarrow\M{\Aa}$ is non degenerate. Therefore one should
think that the quantum space underlying $\Aa$ maps onto a dense subset of the
compact quantum space underlying $\Ap$. This map is a quantum group
homomorphism. This situation is fully analogous to the classical construction
\cite[Sect.~41A--41C]{loom}. One has to keep in mind, however, that in contrast
to the case of classical locally compact groups, our construction is purely
algebraic and thus does not correspond exactly to the classical construction
of Bohr compactification.

Another feature of the classical Bohr compactification $\overline{G}$ of a
locally compact group $G$ is that for any compact group $K$ and a continuous
group homomorphism $\Psi\colon G\to K$ there exists a unique continuous
homomorphism $\overline{\Psi}\colon\overline{G}\to K$ such that
$\Psi=\overline{\Psi}\comp\alpha$ where $\alpha$ is the canonical
homomorphism from $G$ onto a dense subgroup of $\overline{G}$
(cf.~\cite[Sect.~31]{weil}). The analogous
statement is true in our framework. Let $\chi$ denote the inclusion of $\Ap$
into $\M{\Aa}$. Recall that a multiplier Hopf algebra of compact type is
simply a Hopf algebra.

\begin{Thm}\label{univ}
Let $(\Bb,\del_{\Bb})$ be a multiplier Hopf algebra of compact type and let
$\Phi\colon\Bb\to\M{\Aa}$ be a non degenerate homomorphism such that
\begin{equation}\label{hom}
(\Phi\tens\Phi)\comp\del_{\Bb}=\del\comp\Phi
\end{equation}
(where we use the extension of $\del$ to $\M{\Aa}$). Then there exists a
unique Hopf algebra homomorphism $\overline{\Phi}\colon\Bb\to\Ap$ such that
$\Phi=\chi\comp\overline{\Phi}$.
\end{Thm}

\noindent\begin{proof}
Notice first that a non degenerate homomorphism from a unital algebra to
$\M{\Aa}$ must be unital. It follows from \refeq{hom} that the image of $\Phi$
is contained in $\Ap$. Let us define $\overline{\Phi}$ as the same
homomorphism as $\Phi$, but considered now as a map from $\Bb$ to $\Ap$.
This is clearly a Hopf algebra homomorphism and formula
$\Phi=\chi\comp\overline{\Phi}$ is satisfied. The uniqueness of
$\Phi$ follows from the fact that $\chi$ is an embedding.
\end{proof}

The standard reasoning shows that for a discrete quantum group $(\Aa,\del)$ the
Hopf algebra $(\Ap,\del)$ is the unique Hopf algebra with the universal property
described in Theorem \ref{univ}.

Let $(\Aa,\del)$ be a discrete quantum group. In general it is not easy to find
all elements of the algebra $\Ap$, but some of them are easily described.
Recall that an $N$ dimensional corepresentation of $(\Aa,\del)$ is an element
$u\in M_N\tens\M{\Aa}$ such that $(\id\tens\del)u=u_{12}u_{13}$. If
\[
u=\sum_{k,l=1}^{N}e_{k,l}\tens u_{k,l},
\]
where $(e_{k,l})_{k,l=1,\ldots,N}$ are the matrix units in $M_N$, then
\[
\del(u_{k,l})=\sum_{p=1}^{N}u_{k,p}\tens u_{p,l}\in\M{\Aa}\tens\M{\Aa}.
\]
This way we obtain

\begin{Prop}
Let $(\Aa,\del)$ be a discrete quantum group and let $u$ be a finite
dimensional corepresentation of $(\Aa,\del)$ then the matrix elements of $u$
belong to the algebra $\Ap$.
\end{Prop}

It reasonable to conjecture that all almost periodic elements for a discrete
quantum group are linear combinations of matrix elements of finite dimensional
corepresentations. This is the case for classical groups
(cf.~\cite[Sect.~31]{weil}).

It is possible, however, to point out elements that do not belong to $\Ap$. We
shall put this result in the following:

\begin{Prop}\label{niecale}
Let $(\Aa,\del)$ be a discrete quantum group with $\Aa$ infinite dimensional.
Then the algebra $\Ap$ of almost periodic elements does not contain $\Aa$.
\end{Prop}

\noindent\begin{proof}
Let $h$ be the element of $\Aa$ with the property that for any $a\in\Aa$
\[
ah=ha=\eps(a)h.
\]
It is possible to show (cf.~\cite[Sect.~5]{afgd}) that any element $b$ of $\Aa$
has a unique representation in the form
\[
b=(\id\tens\omega)\del(h),
\]
where $\omega\in\Aah$ (notice that for discrete quantum groups
$\Aah\subset\Aa^\st$). However if $\del(h)$ were in $\M{\Aa}\tens\M{\Aa}$ then
the space of right slices of $\del(h)$ with all functionals $\omega\in\Aah$
would have to be finite dimensional.
\end{proof}

If $(\Aa,\del)$ is a multiplier Hopf algebra then the comultiplication,
coinverse and counit extend to the whole algebra $\M{\Aa}$. However in general
$(\M{\Aa},\del)$ is not a Hopf algebra. More precisely we have the following
corollary of Proposition \ref{niecale}.

\begin{Cor}
Let $(\Aa,\del)$ be a discrete quantum group with $\Aa$ infinite dimensional.
Then $(\M{\Aa},\del)$ is not a Hopf algebra.
\end{Cor}

\noindent\begin{proof}
If $(\M{\Aa},\del)$ were a Hopf algebra, in other words, a multiplier Hopf
algebra of compact type, the identity mapping $\M{\Aa}\to\M{\Aa}$ would have
to factor through the inclusion of $\Ap$ into $\M{\Aa}$. Proposition
\ref{niecale} tells us that this is not possible because $\Ap$ does not contain
all elements of $\Aa$ and consequently is not equal to $\M{\Aa}$.
\end{proof}

\section*{acknowledgements}
The author wishes to thank professor Alfons Van Daele for fruitful discussions
and help on the topic of multiplier Hopf algebras and professor Stefaan Vaes
for important comments. He would also like to
express his thanks to professor Joachim Cuntz and colleagues from the
Mathematisches Institut of the University of M\"unster for warm hospitality
and great scientific atmosphere.

\end{document}